\documentclass{amsart}
\usepackage{amssymb}
\usepackage{a4}
\begin{document}
\newtheorem{theorem}{Theorem}[section]
\newtheorem{definition}[theorem]{Definition}
\newtheorem{lemma}[theorem]{Lemma}
\newtheorem{example}[theorem]{Example}
\newtheorem{remark}[theorem]{Remark}
\def\ffrac#1#2{{\textstyle\frac{#1}{#2}}}
\def\rank{\operatorname{Rank}}\def\range{\operatorname{Range}}
\def\oo{\mathfrak{o}}
\def\xx{x^*}\def\XX{X^*}\def\zz{z^*}\def\ZZ{Z^*}\def\tzz{\tilde z^*}\def\tZZ{\tilde Z^*}
\def\ss{s^*}\def\SS{S^*}\def\tss{\tilde s^*}\def\tSS{\tilde S^*}
\def\dss{\ddot s^*}
\makeatletter
 \renewcommand{\theequation}{%
 \thesection.\alph{equation}}
 \@addtoreset{equation}{section}
 \makeatother
\title[Isometry groups of $k$-Curvature homogeneous manifolds]
{Isometry groups of $k$-curvature homogeneous pseudo-Riemannian manifolds}
\author{P. Gilkey and S. Nik\v cevi\'c}
\begin{address}{PG: Mathematics Department, University of Oregon,
Eugene Or 97403 USA.\newline Email: {\it gilkey@darkwing.uoregon.edu}}
\end{address}
\begin{address}{SN: Mathematical Institute, SANU,
Knez Mihailova 35, p.p. 367,
11001 Belgrade,
Serbia and Montenegro.
\newline Email: {\it stanan@mi.sanu.ac.yu}}\end{address}
\begin{abstract} We study the isometry groups of a family of complete $p+2$-curvature homogeneous pseudo-Riemannian
  metrics on $\mathbb{R}^{6+4p}$ which have neutral signature $(3+2p,3+2p)$, and which are $0$-curvature modeled on an
indecomposible symmetric space.\end{abstract}
\keywords{{$k$-curvature homogeneous, homogeneous space, symmetric space, isometry group.}
\newline 2000 {\it Mathematics Subject Classification.} 53B20}
\maketitle
\section{Introduction}\label{sect-1}

Let $\mathcal{M}:=(M,g)$ be a pseudo-Riemannian manifold of signature
$(p,q)$. Let $g_P:=g|_{T_PM}$ (resp.
$\nabla^iR_P:=\nabla^iR|_{T_PM}$) be the restriction of the metric (resp. the $i^{\operatorname{th}}$ covariant derivative of the
curvature tensor) to the tangent space at $P\in M$. We define the {\it $k$-model of $\mathcal{M}$ at $P$} by setting:
$$
\mathfrak{M}_k(\mathcal{M},P):=(T_PM,g_P,R_P,...,\nabla^kR_P)\,.
$$
One says that
$\phi:\mathfrak{M}_k(\mathcal{M}_1,P_1)\rightarrow\mathfrak{M}_k(\mathcal{M}_2,P_2)$ is an {\it isomorphism} from the $k$-model
of $\mathcal{M}_1$ at $P_1$ to the $k$-model of $\mathcal{M}_2$ at $P_2$ if
$\phi$ is a linear isomorphism from $T_{P_1}M_1$ to $T_{P_2}M_2$ with
$$\phi^*g_{2,P_2}=g_{1,P_1}\quad\text{and}\quad
\phi^*\nabla_2^iR_{\mathcal{M}_2,P_2}=\nabla_1^iR_{\mathcal{M}_1,P_1}\quad\text{for}\quad0\le i\le k\,.$$ 
One says that $\mathcal{M}$ is  {\it $k$-curvature homogeneous} if the $k$-models
$\mathfrak{M}_k(\mathcal{M},P)$ and
$\mathfrak{M}_k(\mathcal{M},Q)$ are isomorphic for any $P,Q\in M$.

In the Riemannian setting ($p=0$), Takagi \cite{T74} constructed
$0$-curvature homogeneous complete non-compact Riemannian manifolds; compact examples were exhibited subsequently by Ferus,
Karcher, and M\"unzer
\cite{FKM81}. Although many other examples have been constructed, there are no known
Riemannian manifolds which are $1$-curvature homogeneous but not locally homogeneous and it is natural to conjecture that any
$1$-curvature homogeneous Riemannian manifold is locally homogeneous. 

In the Lorentzian setting ($p=1$), curvature
homogeneous manifolds which are not locally homogeneous were constructed by Cahen et. al.
\cite{CLPT90};
$1$-curvature homogeneous Lorentzian manifolds which are not locally homogeneous have been exhibited by Bueken and Djori\'c
\cite{BD00} and by Bueken and Vanhecke \cite{BV97}. One could conjecture that a $2$-curvature homogeneous Lorentzian manifold must
be locally homogeneous.

It is clear that local homogeneity implies $k$-curvature homogeneity for any
$k$. The following result, due to Singer
\cite{S60} in the Riemannian setting and to F. Podesta and A. Spiro
\cite{PS04} in the general context, provides a partial converse:
\begin{theorem}[Singer, Podesta-Spiro]\label{thm-1.1}
There exists an integer $k_{p,q}$ so that if $\mathcal{M}$ is a complete simply connected pseudo-Riemannian manifold of
signature $(p,q)$ which is $k_{p,q}$-curvature homogeneous, then $(M,g)$ is homogeneous.
\end{theorem}

Sekigawa, Suga, and Vanhecke
\cite{SSV92,SSV95} showed  any $1$-curvature homogeneous complete simply connected Riemannian manifold of dimension $m<5$ is
homogeneous; thus
$k_{0,2}=k_{0,3}=k_{0,4}=1$. The estimate $k_{0,m}<\frac32m-1$ was established by Gromov \cite{Gr86}. Results of \cite{GS04} can
be used to show $k_{p,q}\ge\min(p,q)$; we conjecture $k_{p,q}=\min(p,q)+1$.

If $\mathcal{H}$ is a homogeneous space, let $\mathfrak{M}_k(\mathcal{H}):=\mathfrak{M}_k(\mathcal{H},Q)$ for any point $Q\in
H$; the isomorphism class of $\mathfrak{M}_k(\mathcal{H})$ is independent of the point $Q\in H$. We say that
$\mathcal{M}$ is {\it
$k$-modeled on $\mathcal{H}$} and that $\mathfrak{M}_k(\mathcal{H})$ is {\it a $k$-model for} $\mathcal{M}$ if
$\mathfrak{M}_k(\mathcal{H})$ and $\mathfrak{M}_k(\mathcal{M},P)$ are isomorphic for any $P\in M$.

Throughout this paper, we shall adopt the notational convention that
$$p\ge1\,.$$
 In \cite{GS04x}, we exhibited complete metrics  on
$\mathbb{R}^{6+4p}$ of neutral signature $(3+2p,3+2p)$ which are
$(p+2)$-curvature homogeneous, which are $0$-modeled on an indecomposible symmetric
space, but which are not $(p+3)$-curvature homogeneous; these examples show that the constants $k_{p,q}\rightarrow\infty$ 
as $(p,q)\rightarrow\infty$. The
proof of Theorem
\ref{thm-1.1} rested on a careful analysis of the isometry groups of the model spaces. In this paper, we continue our study of
the manifolds introduced in \cite{GS04x} by examining their isometry groups and the isometry groups of their
$k$-models.

We recall the definition of the metrics on $\mathbb{R}^{6+4p}$ which were introduced in \cite{GS04x}. We will be defining a number
of tensors in this paper and, in the interests of brevity, we shall only give the non-zero components up to the usual symmetries.
Let
$x=(x_1,...,x_m)$ be the usual coordinates on
$\mathbb{R}^m$. Let
$$
  \{x,y,z_1,...,z_p,\tilde y,\tilde z_1,...,\tilde z_p,x^*,y^*,z_1^*,...,z_p^*,
\tilde y^*,\tilde z_1^*,...,\tilde z_p^*\}
$$ 
be coordinates on $\mathbb{R}^{6+4p}$. Let $F=F(y,z_1,...,z_p)\in C^\infty(\mathbb{R}^{p+1})$. Let 
$$\mathcal{M}_{6+4p,F}:=(\mathbb{R}^{6+4p},g_{6+4p,F})$$ where
$g_{6+4p,F}$ is the metric of neutral signature
$(3+2p,3+2p)$ on
$\mathbb{R}^{6+4p}$ with:
\begin{eqnarray*}
&&g_{6+4p,F}(\partial_x,\partial_x)=-2\{F(y,z_1,...,z_p)+y\tilde y+z_1\tilde z_1...+z_p\tilde z_p\},\\
&&g_{6+4p,F}(\partial_x,\partial_{\xx})
  =g_{6+4p,F}(\partial_y,\partial_{y^*})
  =g_{6+4p,F}(\partial_{\tilde y},\partial_{\tilde y^*})=1,\\
&&g_{6+4p,F}(\partial_{z_i},\partial_{\zz_i})
  =g_{6+4p,F}(\partial_{\tilde z_i},\partial_{\tzz_i})=1\,.
\end{eqnarray*}
\begin{theorem}[Gilkey-Nik\v cevi\'c \cite{GS04x}]\label{thm-1.2}
Let $\mathcal{M}=\mathcal{M}_{6+4p,F}$. Then: \begin{enumerate} 
\smallbreak\item All geodesics in $\mathcal{M}$ extend for infinite time.
\smallbreak\item $\exp_P:T_P\mathbb{R}^{6+4p}\rightarrow\mathbb{R}^{6+4p}$ is a
diffeomorphism for all $P\in\mathbb{R}^{6+4p}$.
\smallbreak\item $\nabla^kR(\partial_x,\partial_{\xi_1},\partial_{\xi_2},\partial_x;\partial_{\xi_3},...,\partial_{\xi_{k+2}})
=-\ffrac12(\partial_{\xi_1}\cdot\cdot\cdot\partial_{\xi_{k+2}})g_{6+4p,F}(\partial_x,\partial_x)$
are the non-zero
components of
$\nabla^kR$ where $\xi_i\in\{y,z_1,...,z_p,\tilde y,\tilde z_1,...,\tilde z_p\}$.
\smallbreak\item All scalar Weyl invariants of $\mathcal{M}$ vanish.
\smallbreak\item $\mathcal{M}$ is a symmetric space if and only if $F$ is at most quadratic.
\end{enumerate}
\end{theorem}

\subsection{The manifolds $\mathcal{M}_{6+4p,k}=(\mathbb{R}^{6+4p},g_{6+4p,k})$} We can specialize this construction as follows.
Let $g_{6+4p,k}$ be defined by setting $F=f_{p,k}$ where we let:
\begin{eqnarray*}
&&f_{p,0}(y,z_1,...,z_p):=0,\\
&&f_{p,k}(y,z_1,...,z_p):=z_1y^2+...+z_ky^{k+1}\quad\text{if}\quad1\le k\le p\,.
\end{eqnarray*}
As exceptional cases, we set:
\begin{eqnarray*}
&&f_{p,p+1}(y,z_1,...,z_p):=z_1y^2+...+z_py^{p+1}+y^{p+3},\\
&&f_{p,p+2}(y,z_1,...,z_p):=z_1y^2+...+z_py^{p+1}+e^y\,.
\end{eqnarray*}

\begin{theorem}[Gilkey-Nik\v cevi\'c \cite{GS04x}]\label{thm-1.3}
Let $1\le k\le p+2$.
\ \begin{enumerate}
\smallbreak\item $\mathcal{M}_{6+4p,0}$ is an indecomposible symmetric space.
\smallbreak\item $\mathcal{M}_{6+4p,k}$ is an indecomposible homogeneous space which is not symmetric.
\end{enumerate}
\end{theorem}

\subsection{The manifolds $\mathcal{N}_{6+4p,\psi}=(\mathbb{R}^{6+4p},g_{6+4p,\psi})$} Let $\psi=\psi(y)$ be a real analytic
function of
one variable such that
$$
\psi^{(p+3)}>0,\quad\psi^{(p+4)}>0,\quad\text{and}\quad
\psi^{(p+3)}\ne ae^{by}\,.
$$
Define a metric
$g_{6+4p,\psi}$ on
$\mathbb{R}^{6+4p}$ by taking $F=f_\psi$ where
$$f_\psi(y,z_1,...,z_p):=\psi(y)+z_1y^2+...+z_py^{p+1}\,.$$

The following result shows that the
geometry of a homogeneous pseudo-Riemannian manifold need not determined by the
$k$-model:
\begin{theorem}[Gilkey-Nik\v cevi\'c \cite{GS04x}]\label{thm-1.4}
Let $0\le j<k\le p+2$.
\ \begin{enumerate}
\smallbreak\item $\mathcal{M}_{6+4p,k}$ is $j$-modeled on $\mathcal{M}_{6+4p,j}$;
 $\mathcal{M}_{6+4p,j}$ is not $k$-modeled on $\mathcal{M}_{6+4p,k}$.
\smallbreak\item $\mathcal{N}_{6+4p,\psi}$ is $p+2$-curvature homogeneous and $p+2$-modeled on $\mathcal{M}_{6+4p,p+2}$.
\smallbreak\item $\mathcal{N}_{6+4p,\psi}$ is not $p+3$-curvature homogeneous and not locally homogeneous.
\end{enumerate}
\end{theorem}

\subsection{Isometry groups} Let $G(\mathcal{M})$ (resp. $G(\mathfrak{M}_k)$) be the isometry
group of a pseudo-Riemannian manifold $\mathcal{M}$ (resp. of a $k$-model $\mathfrak{M}_k$). In this paper, we study the groups
$G(\mathcal{M}_{6+4p,k})$, $G(\mathcal{N}_{6+4p,\psi})$, and $G(\mathfrak{M}_k(\mathcal{M}_{6+4p,k},P))$ for any point
$P$ of $\mathbb{R}^{6+4p}$. A byproduct of our study is the following result that shows, not surprisingly, that the symmetric space
$\mathcal{M}_{6+4p,0}$ has the largest isometry group.

\begin{theorem}\label{thm-1.5}
Let $1\le k\le p$. Let $n_p:=(6+4p)+(p+1)(3+2p)+(2p+3)$.
\begin{enumerate}
\smallbreak\item $\dim\{G(\mathcal{M}_{6+4p,0})\}=n_p+(p+1)(2p+1)$.
\smallbreak\item $\dim\{G(\mathcal{M}_{6+4p,k})\}=n_p+(2p+2)+\frac12(2p-k)(2p-k-1)$. 
\smallbreak\item $\dim\{G(\mathcal{M}_{6+4p,p+1})\}=\dim\{G(\mathcal{M}_{6+4p,p})\}-1$.
\smallbreak\item $\dim\{G(\mathcal{M}_{6+4p,p+2})\}=\dim\{G(\mathcal{M}_{6+4p,p+1})\}-1$.
\smallbreak\item $\dim\{G(\mathcal{N}_{6+4p,\psi})\}=\dim\{G(\mathcal{M}_{6+4p,p+2})\}-1$.
\end{enumerate}\end{theorem}

Here is a brief outline to the remainder of this paper. In Section \ref{sect-2}, we review some results from \cite{GS04x}.
In Section \ref{sect-3}, we reduce the proof of Theorem \ref{thm-1.5} to a purely algebraic problem by showing for any
$P\in\mathbb{R}^{6+4p}$ that for $0\le k\le p+2$, we have:
\begin{eqnarray*}
&&\dim\{G(\mathcal{M}_{6+4p,k})\}=6+4p+\dim\{G(\mathfrak{M}_k(\mathcal{M}_{6+4p,k},P))\},\\
&&\dim\{G(\mathcal{N}_{6+4p,\psi})\}=5+4p+\dim\{G(\mathfrak{M}_{p+2}(\mathcal{M}_{6+4p,p+2},P))\}\,.
\end{eqnarray*}
In Section \ref{sect-4}, we complete the proof by determining
$\dim\{G(\mathfrak{M}_k(\mathcal{M}_{6+4p,k},P))\}$ for $0\le k\le p+2$.

\section{Models}\label{sect-2}

It is convenient to work in the purely algebraic setting. Let
$$\mathfrak{M}_\nu:=(V,\langle\cdot,\cdot\rangle,A^0,...,A^\nu)$$
where $\langle\cdot,\cdot\rangle$ is a non-degenerate inner product
of signature $(p,q)$ on a finite dimensional vector space $V$ of dimension $m=p+q$ and where $A^\mu\in\otimes^{4+\mu}V^*$ satisfies
the appropriate symmetries of the covariant derivatives of the curvature tensor for $0\le\mu\le\nu$; if $\nu=\infty$, then the
sequence is infinite. We say that
$\mathfrak{M}_\nu$ is {a \it $\nu$-model} for a pseudo-Riemannian manifold $\mathcal{M}=(M,g)$ if for each point $P\in M$, there
is an isomorphism
$\phi_P:T_PM\rightarrow V$ so that
$$\phi_P^*\langle\cdot,\cdot\rangle=g_P\quad\text{and}\quad \phi_P^*A^\mu=\nabla^\mu R_P\quad\text{for}
\quad 0\le\mu\le \nu\,.$$
Clearly $\mathcal{M}$ is $\nu$-curvature homogeneous if and only if it admits a $\nu$-model.

\subsection{Models for the manifolds $\mathcal{M}_{6+4p,k}$ and $\mathcal{N}_{6+4p,\psi}$}\label{sect-2.1}
Let
$$\mathcal{B}=\{X,Y,Z_1...,Z_p,\tilde Y,\tilde Z_1,...,\tilde Z_p,\XX,Y^*,\ZZ_1,...,\ZZ_p,\tilde Y^*,\tZZ_1,...,\tZZ_p\}$$ 
be a basis for $\mathbb{R}^{6+4p}$. Define a hyperbolic inner-product on $\mathbb{R}^{6+4p}$ by pairing ordinary variables with
the corresponding dual $\star$-variables:
\begin{equation}\label{eqn-2.a}
\langle X,\XX \rangle=
\langle Y,Y^*\rangle=\langle\tilde Y,\tilde Y^*\rangle=
\langle{Z_i},{\ZZ_i}\rangle=\langle{\tilde Z_i},{\tZZ_i}\rangle=1\,.
\end{equation}
Define $A^0\in\otimes^4(\mathbb{R}^{6+4p})^*$ with non-zero components:
$$A^0(X,Y,\tilde Y,X)=A^0(X,{Z_i},{\tilde Z_i},X)=1\,.
$$
Define tensors $A^i\in\otimes^{4+i}(\mathbb{R}^{6+4p})^*$ for $1\le i\le p$ with non-zero components:
\begin{eqnarray*}
&&A^i(X,Y,{Z_i},X;Y,...,Y)=1,\\
&&A^i(X,Y,Y,X;{Z_i},Y,...,Y)=1,...,\\
&&A^i(X,Y,Y,X;Y,...,Y,{Z_i})=1\,.
\end{eqnarray*}
Finally define $A^{p+1}\in\otimes^{5+p}(\mathbb{R}^{6+4p})^*$ and $A^{p+2}\in\otimes^{6+p}(\mathbb{R}^{6+4p})^*$ by setting
\begin{eqnarray*}
&&A^{p+1}(X,Y,Y,X;Y,...,Y)=1,\\
&&A^{p+2}(X,Y,Y,X;Y,...,Y)=1\,.
\end{eqnarray*}
Define models:
$$\mathfrak{M}_{6+4p,k}:=(\mathbb{R}^{6+4p},\langle\cdot,\cdot\rangle,A^0,...,A^k)\quad\text{for}\quad 0\le k\le p+2\,.$$

\begin{lemma}[Gilkey-Nik\v cevi\'c \cite{GS04x}]\label{lem-2.1}Let $0\le k\le p+2$. \begin{enumerate}
\smallbreak\item $\mathfrak{M}_{6+4p,k}$ is a $k$-model for $\mathcal{M}_{6+4p,k}$.
\smallbreak\item $\mathfrak{M}_{6+4p,{p+2}}$ is a $p+2$-model for $\mathcal{N}_{6+4p,\psi}$.
\end{enumerate}
\end{lemma} 
 
\section{Isometry groups in the geometric setting}\label{sect-3}
In this section we will reduce the proof of Theorem \ref{thm-1.5} to a purely algebraic problem by showing:

\begin{theorem}\label{thm-3.1}Let $0\le k\le p+2$. \begin{enumerate}
\smallbreak\item $\dim\{G(\mathcal{M}_{6+4p,k})\}=6+4p+\dim\{G(\mathfrak{M}_{6+4p,k})\}$.
\smallbreak\item $\dim\{G(\mathcal{N}_{6+4p,\psi})\}=5+4p+\dim\{G(\mathfrak{M}_{6+4p,p+2})\}$.
\end{enumerate}
\end{theorem}

The proof of Theorem \ref{thm-3.1} will be based on several Lemmas. In Lemma \ref{lem-3.2}, we review a basic result about group
actions. In Lemma
\ref{lem-3.3}, we relate the full isometry group
$G(\cdot)$ to the isotropy subgroup. In Lemma \ref{lem-3.4}, we relate the isotropy subgroup to the isometry group of the
$\infty$-model. In Lemma
\ref{lem-3.5}, we relate
isometry group of the $\infty$-model to the isometry group of an appropriate finite model.

The following result is well known.

\begin{lemma}\label{lem-3.2}
Let $G$ be a Lie group which acts continuously on a metric space $X$. If $x\in X$, let $G\cdot x$ be
the orbit and let
$G_x=\{g\in G:gx=x\}$ be the isotropy subgroup.\begin{enumerate}
\smallbreak\item We have a smooth principle bundle $G_x\rightarrow
G\rightarrow G\cdot x$.
\smallbreak\item $\dim\{G\}=\dim\{G_x\}+\dim\{G\cdot x\}$.
\end{enumerate}\end{lemma}

We can relate $\dim\{G(\mathcal{M})\}$ to $\dim\{G_P(\mathcal{M})\}$ for $\mathcal{M}=\mathcal{M}_{6+4p,k}$ or
$\mathcal{M}=\mathcal{N}_{6+4p,\psi}$.

\begin{lemma}\label{lem-3.3}
Let $P\in\mathbb{R}^{6+4p}$. Let $0\le k\le p+2$.
\ \begin{enumerate}
\smallbreak\item $\dim\{G(\mathcal{M}_{6+4p,k})\}=6+4p+\dim\{G_P(\mathcal{M}_{6+4p,k})\}$.
\smallbreak\item $\dim\{G(\mathcal{N}_{6+4p,\psi})\}=6+4p-1+\dim\{G_P(\mathcal{N}_{6+4p,\psi})\}$.
\end{enumerate}\end{lemma}

\begin{proof} We apply Lemma \ref{lem-3.2} to the canonical action of $G(\mathcal{M})$ on $\mathbb{R}^{6+4p}$. Assertion (1)
follows as
$\mathcal{M}_{6+4p,k}$ is a homogeneous space. Let $\nu\ge2$. Set
$$\alpha_{6+4p,\nu}(\psi):=\psi^{(\nu+p+3)}\{\psi^{(p+3)}\}^{\nu-1}\{\psi^{(p+4)}\}^{-\nu}\,.$$
We showed \cite{GS04x} that if $\mathcal{B}$ is a basis satisfying the normalizations of Section \ref{sect-2.1}, then
the only non-zero components of $\nabla^{\nu+p+1}R$ are given by:
\begin{equation}\label{eqn-3.a}
\nabla^{\nu+p+1}R(X,Y,Y,X;Y,...,Y)=\alpha_{6+4p,\nu}(\psi)\,.
\end{equation}
We also showed  that the following assertions are equivalent:
\begin{enumerate}
\smallbreak\item $\alpha_{6+4p,\nu}(\psi_1)(P_1)=\alpha_{6+4p,\nu}(\psi_2)(P_2)$ for all $\nu\ge2$.
\smallbreak\item There exists an isometry
$\phi:\mathcal{N}_{6+4p,\psi_1}\rightarrow\mathcal{N}_{6+4p,\psi_2}$ with $\phi(P_1)=P_2$.
\end{enumerate}
The functions $\alpha_{6+4p,\nu}(\psi)$ are constant on the hyperplanes $y=c$; thus the group of isometries acts transitively on
such a hyperplane. Consequently 
$$\dim\{G(\mathcal{N}_{6+4p,\psi})\}\ge\dim\{G_P(\mathcal{N}_{6+4p,\psi})\}+6+4p-1\,.$$
Since $\mathcal{N}_{6+4p,\psi}$ is not a homogeneous space, equality holds.
\end{proof}

Let $P\in M$. We can show that $G_P(\mathcal{M})$ is isomorphic to $G(\mathfrak{M}_\infty(\mathcal{M},P))$ under certain
circumstances.

\begin{lemma}\label{lem-3.4}\ \begin{enumerate}\item Let $\mathcal{M}_1:=(M_1,g_1)$ and $\mathcal{M}_2:=(M_2,g_2)$ be real
analytic. Assume for $\varrho=1,2$ that there are points $P_\varrho\in M_\varrho$ so
$\exp_{P_\varrho}:T_{P_\varrho}M_\varrho\rightarrow M_\varrho$ is a diffeomorphism. If $\phi:
T_{P_1}M_1\rightarrow T_{P_2}M_2$ induces an isomorphism from $\mathfrak{M}_\infty(\mathcal{M}_1,P_1)$ to
$\mathfrak{M}_\infty(\mathcal{M}_2,P_2)$, then
$\Phi:=\exp_{P_2}\circ\phi\circ\exp_{P_1}^{-1}$ is an isometry from $\mathcal{M}_1$ to $\mathcal{M}_2$.
\smallbreak\item If $\mathcal{M}=\mathcal{M}_{6+4p,k}$ or if
$\mathcal{M}=\mathcal{N}_{6+4p,\psi}$, then $G_{P}(\mathcal{M})=G(\mathfrak{M}_\infty(\mathcal{M},P))$ for any point
$P\in\mathbb{R}^{6+4p}$.
\end{enumerate}
\end{lemma}

\begin{proof} Belger and Kowalski \cite{BeKo94} note about analytic pseudo-Riemannian metrics that the
``metric
$g$ is uniquely determined, up to local isometry, by the tensors $R$, $\nabla R$, ..., $\nabla^kR$, ... at one point.''; see also
Gray \cite{Gr73} for related work. The first assertion now follows; the second follows immediately from the first and from Theorem
\ref{thm-1.2}.
\end{proof}

We now replace the infinite model by a finite model:

\begin{lemma}\label{lem-3.5}
Let $P\in\mathbb{R}^{6+4p}$. Let $0\le k\le p+2$. Then:
\begin{enumerate}
\smallbreak\item $G(\mathfrak{M}_\infty(\mathcal{M}_{6+4p,k},P))=G(\mathfrak{M}_{6+4p,k})$.
\smallbreak\item $G(\mathfrak{M}_\infty(\mathcal{N}_{6+4p,\psi},P))=G(\mathfrak{M}_{6+4p,p+2})$.
\end{enumerate}
\end{lemma} 

\begin{proof} If $\mathcal{M}$ is a pseudo-Riemannian manifold, restriction induces an injective map
$$
  r:G(\mathfrak{M}_\infty(\mathcal{M},P))\rightarrow G(\mathfrak{M}_k(\mathcal{M},P))\,.
$$
Suppose that $\mathcal{M}=\mathcal{M}_{4p+6,k}$ for $k<p+2$. Then $\nabla^jR=0$ for $j>k$; consequently any isomorphism of the
$k$-model is an isomorphism of the $\infty$-model; this proves Assertion (1) for $0\le k\le p+1$.

To deal with the remaining cases, we suppose that
$\psi^{(p+3)}$ and $\psi^{(p+4)}$ are always positive, but drop the restriction that $\psi^{(p+3)}\ne ae^{by}$. Choose a basis
$\mathcal{B}$ for $T_PM$ satisfying the normalizations of Section \ref{sect-2.1}. If
$g\in G(\mathfrak{M}_{p+2}(\mathcal{M}_{6+4p,p+2},P))$, then
$g\mathcal{B}$ also satisfies the normalizations of Section \ref{sect-2.1}. We may then apply Equation (\ref{eqn-3.a}) to see that
$g$ is in fact an isomorphism of the $\infty$-model since $g$ preserves $\nabla^kR$ for any $k>p+2$. The first assertion with
$k=p+2$ and the second assertion of the Lemma now follow; this also completes the proof of Theorem \ref{thm-3.1}.
\end{proof}

\section{Isometry groups of the models}\label{sect-4}
Let $\mathbb{R}^{3+2p}:=\operatorname{Span}\{X,Y,Z_1,...,Z_p,\tilde Y\tilde Z_1,...,\tilde Z_p\}$ and let
$B^i\in\otimes^{4+i}(\mathbb{R}^{3+2p})^*$ be the restriction of $A^i$ to $\mathbb{R}^{3+2p}$. We introduce
the affine models by restricting the domain and suppressing the metric:
$$\mathfrak{A}_{3+2p,k}:=(\mathbb{R}^{3+2p},B^0,...,B^k)\,.$$

\begin{lemma}\label{lem-4.1}
$\dim\{G(\mathfrak{M}_{6+4p,k})\}=\dim\{G(\mathfrak{A}_{3+2p,k})\}+(p+1)(3+2p)$.
\end{lemma}

\begin{proof} Let $\oo(s)$ be Lie algebra of skew-symmetric $s\times s$ real matrices. Set
\begin{eqnarray*}
\mathcal{S}:&=&(S_1,...,S_{3+2p})=(X,Y,Z_1...,Z_p,\tilde Y,\tilde Z_1...,\tilde Z_p),\\
\mathcal{S}^*:&=&(\SS_1,...,\SS_{3+2p})=(\XX,Y^*,\ZZ_1,...,\ZZ_p,\tilde Y^*,\tZZ_1,...,\tZZ_p),\\
\mathcal{K}:&=&\{\xi\in\mathbb{R}^{6+4p}:A^0(\xi,\eta_1,\eta_2,\eta_3)=0
\ \forall\ \eta_i\in\mathbb{R}^{6+4p}\}\\
&=&\operatorname{Span}\{\SS_1,...,\SS_{3+2p}\}\,.
\end{eqnarray*}

Let $g\in G(\mathfrak{M}_{6+4p,k})$. The space $\mathcal{K}$ is preserved by $g$. Thus
$$gS_i=\textstyle\sum_{i,j}\{g_{0,ij}S_j+g_{1,ij}\SS_j\}\quad\text{and}\quad
g\SS_i=\textstyle\sum_{i,j}\{g_{2,ij}\SS_j\}\,.
$$
By Equation (\ref{eqn-2.a}), $\langle gS_i,gS_j\rangle=0$ and $\langle gS_i,g\SS_j\rangle=\delta_{ij}$. Thus
$$
\textstyle\sum_k\{g_{0,ik}g_{1,jk}+g_{1,ik}g_{0,jk}\}=0\quad\text{and}\quad
\textstyle\sum_k\{g_{0,ik}g_{2,jk}\}=\delta_{ij}\,.
$$
for all $i,j$. Set $\gamma:=g_0g_1^t$. One then has
\begin{equation}\label{eqn-4.a}
g_0\in G(\mathfrak{A}_{3+2p,k}),\quad \gamma+\gamma^t=0,\quad\text{and}\quad g_0g_2^t=\operatorname{id}\,.
\end{equation}
Conversely, if Equation (\ref{eqn-4.a}) is satisfied then $g\in G(\mathfrak{M}_{6+4p,k})$.
The map $g\rightarrow(g_0,\gamma)$ yields an identification of
$$G(\mathfrak{M}_{6+4p,k})=G(\mathfrak{A}_{3+2p,k})\times\oo(3+2p)$$
as a twisted product. The Lemma follows as
$\dim\{\oo(3+2p)\}=\frac12(3+2p)(2+2p)$.
\end{proof}

There is a natural action of $G(\mathfrak{A}_{3+2p,k})$ on $\mathbb{R}^{3+2p}$. We continue our study by relating
$G(\mathfrak{A}_{3+2p,k})$ and the isotropy subgroup
$G_X(\mathfrak{A}_{3+2p,k})$.

\begin{lemma}\label{lem-4.2}
\ \begin{enumerate}
\smallbreak\item 
 $\dim\{G(\mathfrak{A}_{3+2p,k})\}=\dim\{G_X(\mathfrak{A}_{3+2p,k})\}+2p+3$ for $k\le p+1$.
\smallbreak\item 
$\dim\{G(\mathfrak{A}_{3+2p,p+2})\}=\dim\{G_X(\mathfrak{A}_{3+2p,p+2})\}+2p+2$.
\end{enumerate}
\end{lemma}

\begin{proof} Lemma \ref{lem-4.2} will follow from Lemma \ref{lem-3.2} and the following relations:
\begin{equation}\label{eqn-4.b}
\begin{array}{l}
G(\mathfrak{A}_{3+2p,k})X=\{\xi\in\mathbb{R}^{3+2p}:\langle\xi,X^*\rangle\ne0\}\text{\ \  if\ \  }k\le p+1,\\
G(\mathfrak{A}_{3+2p,p+2})X=\{\xi\in\mathbb{R}^{3+2p}:\langle\xi,X^*\rangle=\pm1\}\,.\vphantom{\vrule height 11pt}
\end{array}\end{equation}

We first show $\supset$ holds in Equation (\ref{eqn-4.b}). Let $\xi\in\mathbb{R}^{3+2p}$. Assume that
$$a:=\langle\xi,\XX\rangle\ne0\,.$$
Set $gX=\xi$ and set
$$\begin{array}{lll}
\varepsilon_0:=(a^2)^{-1/(p+3)},&gY:=\varepsilon_0Y,&g\tilde Y:=a^{-2}\varepsilon_0^{-1}\tilde Y,\\ 
\varepsilon_i:=\{a^2\varepsilon_0^{i+1}\}^{-1},&gZ_i:=\varepsilon_iZ_i,&g\ZZ_i:=\varepsilon_i^{-1}a^{-2}\tilde Z_i\,.
\end{array}$$
The non-zero components of $\nabla^iR$ for $1\le i\le p+2$ are then given by
\begin{eqnarray*}
&&R(gX,gY,g\tilde Y,gX)=a^2\varepsilon_0a^{-2}\varepsilon_0^{-1}=1,\\
&&R(gX,gZ_i,g\tilde Z_i,gX)=a^2\varepsilon_i\varepsilon_i^{-1}a^{-2}=1,\\
&&\nabla R(gX,gY,gZ_1,gX;gY)=\nabla R(gX,gY,gY,gX;gZ_1)=a^2\varepsilon_0^2\varepsilon_1=1,...\\
&&\nabla^pR(gX,gY,gZ_p,gX;gY,...,gY)=\nabla^pR(gX,gY,gY,gX;gZ_p,gY,...,gY)=...\\
&&\quad=\nabla^pR(gX,gY,gY,gX;gY,...,gY,gZ_p)=a^2\varepsilon_0^{p+1}\varepsilon_p=1,\\
&&\nabla^{p+1}R(gX,gY,gY,gX;gY,...,gY)=a^2\varepsilon_0^{p+3}=1,\\
&&\nabla^{p+2}R(gX,gY,gY,gX;gY,...,gY)=a^2\varepsilon_0^{p+4}=\varepsilon_0\,.
\end{eqnarray*}
Thus $g\in G(\mathfrak{A}_{3+2p,p+1})$. Furthermore, $g\in G(\mathfrak{A}_{3+2p,p+2})$ if $a^2=1$. Consequently:
\begin{equation}\label{eqn-4.c}
\begin{array}{l}
\{\xi\in\mathbb{R}^{3+2p}:\langle\xi,X^*\rangle\ne0\}\subset G(\mathfrak{A}_{3+2p,k})\cdot X\quad\text{for}\quad k\le p+1,\\
\{\xi\in\mathbb{R}^{3+2p}:\langle\xi,X^*\rangle=\pm1\}\subset G(\mathfrak{A}_{3+2p,p+2})\cdot X\,.\vphantom{\vrule height 12pt}
\end{array}\end{equation}

We must establish the reverse inclusions to complete the proof.
Let $\xi\in \mathbb{R}^{3+2p}$. Let
$J_\xi(\eta_1,\eta_2):=R(\xi,\eta_1,\eta_2,\xi)$ be the {\it Jacobi form}.
Adopt the Einstein convention and sum over repeated indices to expand
$$\textstyle\xi=aX+b^iZ_i+\tilde b^i\tilde Z_i$$
 where $a=\langle\xi,\XX\rangle$. We have the following cases
\begin{enumerate}
\item If $a=0$, then $J_\xi=0$ on $\operatorname{Span}\{Y,\tilde Y,Z_i,\tilde Z_i\}$ so
  $\rank(J_\xi)\le 1$.
\item If $a\ne0$, then $J_\xi(Y,\tilde Y)\ne0$ so $\rank(J_\xi)\ge 2$.
\end{enumerate}
If $g\in G(\mathfrak{A}_{3+2p,k})$, then $\rank\{J_\xi\}=\rank\{J_{g\xi}\}$. Consequently
$$
\langle\xi,\XX\rangle=0\Leftrightarrow\rank(J_\xi)\le 1\Leftrightarrow\rank(J_{g\xi})\le 1\Leftrightarrow\langle
g\xi,\XX\rangle=0
$$
Consequently we have
\begin{equation}\label{eqn-4.d}
\begin{array}{l}
G(\mathfrak{A}_{3+2p,k})\cdot X\subset\{\xi\in\mathbb{R}^{3+2p}:\langle\xi,X^*\rangle\ne0\},\\
G(\mathfrak{A}_{3+2p,k})\cdot \operatorname{Span}\{Y,Z_i,\tilde Z_i\}=\operatorname{Span}\{Y,Z_i,\tilde Z_i\}\,.
\vphantom{\vrule height 11pt}
\end{array}\end{equation}
Suppose $k=p+2$. Since $\rank(J_Y)=0$, $\rank(J_{gY})=0$ so $\langle gY,X^*\rangle=0$. Expand
\begin{eqnarray*}
&&gX=aX+a_0Y+\tilde a_0\tilde Y+a^iZ_i+\tilde a^i\tilde Z_i,\\
&&gY=\phantom{aX+.}b^0Y+\tilde b^0\tilde Y+b^iZ_i+\tilde b^i\tilde Z_i\,.
\end{eqnarray*}
Then
\begin{eqnarray*}
&&1=\nabla^{p+1}R(gX,gY,gY,gX;gY,...,gY)=a^2(b^0)^{p+3},\\
&&1=\nabla^{p+2}R(gX,gY,gY,gX;gY,...,gY)=a^2(b^0)^{p+4}\,.
\end{eqnarray*}
This shows that $a^2=1$ and $b^0=1$ so
\begin{equation}\label{eqn-4.e}
\begin{array}{l}
G(\mathfrak{A}_{3+2p,p+2})X\subset\{\xi\in\mathbb{R}^{3+2p}:\langle\xi,X^*\rangle=\pm1\},\\
G(\mathfrak{A}_{3+2p,p+2})Y\subset\{\xi\in\mathbb{R}^{3+2p}:\langle\xi,X^*\rangle=0,\text{ and }\langle\xi,Y^*\rangle=1\}\,.
\vphantom{\vrule height 11pt}
\end{array}\end{equation}
Equations (\ref{eqn-4.c}), (\ref{eqn-4.d}), and (\ref{eqn-4.e}) now imply Equation (\ref{eqn-4.b}); the
Lemma follows.
\end{proof}

We now consider the double isotropy group
$$G_{X,Y}(\mathfrak{A}_{3+2p,k})=\{g\in G(\mathfrak{A}_{3+2p,k}):gX=X\text{ and }gY=Y\}\,.$$

\begin{lemma}\label{lem-4.3}
\ \begin{enumerate}
\smallbreak\item $\dim\{G_X(\mathfrak{A}_{3+2p,0})\}=(p+1)(2p+1)$.
\smallbreak\item 
$\dim\{G_X(\mathfrak{A}_{3+2p,k})\}=\dim\{G_{X,Y}(\mathfrak{A}_{3+2p,k})\}+2p+2$ for $1\le k\le p$.
\smallbreak\item $\dim\{G_X(\mathfrak{A}_{3+2p,k})\}=\dim\{G_{X,Y}(\mathfrak{A}_{3+2p,k})\}+2p+1$ for $k=p+1,p+2$.
\smallbreak\item $G_{X,Y}(\mathfrak{A}_{3+2p,p})=G_{X,Y}(\mathfrak{A}_{3+2p,p+1})=G_{X,Y}(\mathfrak{A}_{3+2p,p+2})$.
\end{enumerate}\end{lemma}

\begin{proof} As noted above, the Jacobi form $J_X(\cdot,\cdot)=R(X,\cdot,\cdot,X)$ defines a
non-singular bilinear form of signature $(p+1,p+1)$ on
$$W:=\operatorname{Span}\{Y,Z_1,...,Z_p,\tilde Y,\tilde Z_1,...,\tilde Z_p\}=\{\xi:\operatorname{Rank}(J_\xi)\le 1\}\,.$$
Let $O(W,J_X)$ be the associated orthogonal group. If $g\in G_X(\mathfrak{A}_{3+2p,k})$, then we have $gW=W$ by Equation
(\ref{eqn-4.d}). Since
$gX=X$, we may safely identify $g$ with $g|_W$. Furthermore,
$$J_X(\xi,\eta)=J_{gX}(g\xi,g\eta)=J_X(g\xi,g\eta)\quad\text{so}\quad
G_X(\mathfrak{A}_{3+2p,k})\subset O(W,J_X)\,.$$
Conversely, if $g$ is a linear map of $W$ which preserves $J_X$, we may extend $g$ to $\mathbb{R}^{3+2p}$ by defining $gX=X$ and
thereby obtain an element of
$G_X(\mathfrak{A}_{3+2p,0})$. Thus $G_X(\mathfrak{A}_{3+2p,0})=O(W,J_X)$.
Assertion (1) now follows since
$$\dim\{O(W,J_X)\}=\ffrac12\dim W(\dim W-1)=\ffrac12(1+2p)(2+2p)\,.$$

Assertions (2) and (3) will follow from Lemma \ref{lem-3.2} and from the relations:
\begin{equation}\label{eqn-4.f}
\begin{array}{l}
G_{X}(\mathfrak{A}_{3+2p,k})\cdot Y=\{\xi\in W:\langle\xi,Y^*\rangle\ne0\}\quad\text{for}\quad1\le k\le p,\\
G_{X}(\mathfrak{A}_{3+2p,p+1})\cdot Y=\{\xi\in W:\langle\xi,Y^*\rangle^{p+3}=1\},\vphantom{\vrule height 12pt}\\
G_{X}(\mathfrak{A}_{3+2p,p+2})\cdot Y=\{\xi\in W:\langle\xi,Y^*\rangle=1\}\,.\vphantom{\vrule height 12pt}
\end{array}\end{equation}

If $\xi\in W$, let $S_\xi(\eta):=\nabla R(X,\xi,\xi,X;\eta)$. Expand
\begin{equation}\label{eqn-4.g}
\xi=b^0Y+\tilde b^0\tilde Y+b^iZ_i+\tilde b^i\tilde Z_i\,.
\end{equation}
We then have that
\begin{eqnarray*}
&&S_\xi(X)=0,\quad S_\xi(\tilde Z_i)=0,\quad S_\xi(Y)=2b^0b^1,\\
&& S_\xi(Z_1)=(b^0)^2,\quad\text{and}\quad
S_\xi(Z_i)=0\quad\text{for}\quad i\ge2\,.
\end{eqnarray*}
Thus $S_\xi=0$ if and only if $b^0=\langle\xi,Y^*\rangle=0$. It now follows that for $k\ge1$ we have
\begin{equation}\label{eqn-4.h}
\begin{array}{l}
G_X(\mathfrak{A}_{3+2p,k})Y\subset\{\xi\in W:\langle\xi,Y^*\rangle\ne0\},\\
G_X(\mathfrak{A}_{3+2p,k})\operatorname{Span}\{Z_i,\tilde Y,\tilde Z_i\}
  \subset\operatorname{Span}\{Z_i,\tilde Y,\tilde Z_i\}\,.
\vphantom{\vrule height 11pt}\end{array}
\end{equation}
Since $a=1$, the analysis used to prove Lemma \ref{lem-4.2} shows $(b^0)^{p+3}=1$ if $k=p+1$ and $b^0=1$ if $k=p+2$. This
establishes the inclusions $\subset$ in Equation (\ref{eqn-4.f}).

We complete the proof by establishing the reverse inclusions in Equation (\ref{eqn-4.f}). Expand $\xi$ in the form given in
Equation (\ref{eqn-4.g}). Assume
$b^0\ne0$. Let
$gX=X$,
$gY=\xi$,
$g\tilde Y=(b^0)^{-1}\tilde Y$,
$$gZ_i:=\varepsilon_i\{Z_i-(b^0)^{-1}\tilde b^i\tilde Y\}\quad\text{and}\quad
  g\tilde Z_i:=\varepsilon_i^{-1}\{\tilde Z_i-(b^0)^{-1}b^i\tilde Y\}\,.
$$
The possibly non-zero components of $R$ are then given by
\begin{eqnarray*}
&&R(gX,gY,g\tilde Y ,gX)=1,\\
&&R(gX,gY,gZ_i,gX)=\varepsilon_i\{\tilde b^i-(b^0)(b^0)^{-1}\tilde b^i\}=0,\\
&&R(gX,gY,g\tilde Z_i,gX)=\varepsilon_i^{-1}\{b^i-(b^0)(b^0)^{-1}b^i\}=0,\\
&&R(gX,gZ_i,g\tilde Z_i,gX)=\varepsilon_i^{-1}\varepsilon_i=1\,.
\end{eqnarray*}
The non-zero components of $\nabla^iR$ for $1\le i\le p$ are given by
\begin{eqnarray*}
&&\nabla^iR(gX,gY,gZ_i,gX;gY,...,gY)=...\\
&=&\nabla^iR(gX,gY,gY,gX;gY,...,gZ_i)=(b^0)^{i+1}\varepsilon_i\,.
\end{eqnarray*}
We therefore set $\varepsilon_i=(b^0)^{-i-1}$ for $1\le i\le p$ to ensure $g\in G(\mathfrak{A}_{3+2p,p})$. 

The non-zero components of
$\nabla^iR$ for
$i=p+1,p+2$ are
$$
\nabla^iR(gX,gY,gY,gX;gY,...,gY)=(b^0)^{i+2}\,.
$$
If $(b^0)^{p+3}=1$, then $g\in G(\mathfrak{A}_{3+2p,p+1})$; if $b^0=1$, then $g\in G(\mathfrak{A}_{3+2p,p+2})$. This establishes
the reverse inclusions in Equation (\ref{eqn-4.f}) and completes the proof of Assertions (2) and (3); Assertion (4) is immediate.
\end{proof}

Let $W(p):=\operatorname{Span}\{Z_1,...,Z_p,\tilde Z_1,...,\tilde Z_p\}$. Let 
$\{\beta_1,...,\beta_p,\tilde\beta_1,...,\tilde\beta_p\}$ be the
corresponding dual basis for the dual space $\mathcal{W}(p):=W(p)^*$. The curvature tensor $R(X,\cdot,\cdot,X)$ defines a
non-degenerate form $\langle\cdot,\cdot\rangle$ on
$W(p)$; dually on $\mathcal{W}(p)$ we have:
$$\langle\beta_i,\beta_j\rangle=\langle\tilde\beta_i,\tilde\beta_j\rangle=0,
  \quad\langle\beta_i,\tilde\beta_j\rangle=\delta_{ij}\,.$$
Let
$\mathcal{O}(p)$ be the associated orthogonal group on $\mathcal{W}(p)$. Let
$$\mathcal{O}(p,k):=\{h\in \mathcal{O}(p):h\beta_i=\beta_i\quad\text{for}\quad 1\le i\le k\}$$
be the simultaneous isotropy group. We set $\mathcal{O}(p,0)=\mathcal{O}(p)$. Theorem
\ref{thm-1.5} will now follow from the following result:

\begin{lemma}\label{lem-4.4} Let $1\le k\le p$.
\begin{enumerate}
\item $G_{X,Y}(\mathfrak{A}_{3+2p,k})=\mathcal{O}(p,k)$.
\item $\mathcal{O}_{\tilde\beta_1}(p,k)=\mathcal{O}(p-1,k-1)$.
\item $\dim\{\mathcal{O}(p,k)\}=\dim\{\mathcal{O}(p-1,k-1)\}+2p-k-1$.
\item $\dim\{\mathcal{O}(p,k)\}=\frac12(2p-k)(2p-k-1)$.
\end{enumerate}\end{lemma}

\begin{proof} Let $g\in G_{X,Y}(\mathfrak{A}_{3+2p,k})$. Let $\xi\in\operatorname{Span}\{Z_1,...,Z_p,\tilde Y,\tilde
Z_1,...,\tilde Z_p\}$. We may use Equation (\ref{eqn-4.h}) and the relation
$R(X,Y,g\xi,X)=R(X,Y,\xi,X)$, to see
$$g\tilde Y=\tilde Y+a^iZ_i+a^{\tilde i}\tilde Z_i,\quad
  gZ_i=a_i^jZ_j+a_i^{\tilde j}\tilde Z_{\tilde j},\quad
  g\tilde Z_{\tilde i}=a_{\tilde i}^jZ_j+a_{\tilde i}^{\tilde j}\tilde Z_{\tilde j}\,.
$$
Consequently $\operatorname{Span}_{1\le i\le p}\{gZ_i,g\tilde Z_{\tilde i}\}=\operatorname{Span}_{1\le i\le p}\{Z_i,\tilde
Z_{\tilde i}\}$ and the relation
$$R(X,gZ_i,g\tilde Y,X)=R(X,g\tilde Z_{\tilde i},g\tilde Y,X)=0$$ 
implies $a^i=a^{\tilde i}=0$. Thus $g\tilde Y=\tilde Y$ and $g:W(p)\rightarrow W(p)$; this shows that $g$ is determined by
its restriction to ${W(p)}$. Let $h:={}^*g$ denote the dual action of $g$ on $\mathcal{W}(p)$. The isomorphism of Assertion (1) now
follows as:
\begin{eqnarray*}
&&R(X,g\xi_1,g\xi_2,R)=R(X,\xi_1,\xi_2,X)\ \forall\xi_1,\xi_2\Leftrightarrow h\in \mathcal{O}(p)\,,\\
&&\nabla^iR(X,Y,g\xi,X;Y,...,Y)=\nabla^iR(X,Y,\xi,X;Y,...,Y)\ \forall\xi\Leftrightarrow h\beta_i=\beta_i\,.
\end{eqnarray*} 

If $h(\beta_1)=\beta_1$ and $h(\tilde\beta_1)=\tilde\beta_1$, then $h$ preserves 
$$\operatorname{Span}\{\beta_1,\tilde\beta_1\}^\perp=
\operatorname{Span}\{\beta_2,...,\beta_p,\tilde\beta_2,...,\tilde\beta_p\}\,.$$
The isomorphism of Assertion (2) now follows by restricting $h$ to this subspace and by renumbering the variables appropriately.

We set
$$\mathcal{W}(p,k):=\{\xi\in\mathcal{W}(p):\langle\xi,\xi\rangle=0,\ 
    \langle\xi,\beta_1\rangle=1,\ 
    \langle\xi,\beta_i\rangle=0\text{ for }2\le i\le k\}\,.$$
If $h\in\mathcal{O}(p,k)$, then $h$ preserves $\langle\cdot,\cdot\rangle$ and $h$ preserves $\{\beta_1,...,\beta_k\}$.
Consequently $h\tilde\beta_1\in\mathcal{W}(p,k)$ as $\tilde\beta_1$ satisfies these relations. Conversely,
$\xi\in\mathcal{W}(p,k)$ if and only if
$$\xi=b^1\beta_1+\sum_{1<i}b^i\beta_i+\tilde\beta_1+\sum_{k<i}\tilde b^i\tilde\beta_i\quad\text{where}\quad
    b^1+\sum_{k<i}b^i\tilde b^i=0\,.$$
Since the variables $\{b^2,...,b^p,\tilde b^{k+1},...,\tilde b^p\}$ can be chosen arbitrarily,
$$\mathcal{W}(p,k)=\mathbb{R}^{p-1+p-k}\quad\text{so}\quad\dim\mathcal{W}(p,k)=2p-k-1\,.$$
We show that $\xi\in\mathcal{O}(p,k)\tilde\beta_1$ by finding $h\in\mathcal{O}(p,k)$ so $h\tilde\beta_1=\xi$. Set:
$$\begin{array}{llll}
h\beta_i=\beta_i&\text{for }1\le i\le k,&
h\beta_i=\beta_i-\tilde b^i\beta_1&\text{for }k<i,\\
h\tilde\beta_1=\xi,&&
h\tilde\beta_i=\tilde\beta_i-b^i\beta_1&\text{for }1<i\,.\vphantom{\vrule height 11pt}
\end{array}$$
This shows $\mathcal{O}(p,k)\cdot\tilde\beta_1=\mathcal{W}(p,k)$. Assertion (3) now follows from Assertion (2) and from Lemma
\ref{lem-3.2}.

As $\dim\{\mathcal{O}(p-k)\}=\frac12(2p-2k)(2p-2k-1)$, Assertion (4) follows by induction. \end{proof}

\section*{Acknowledgments} Research of P. Gilkey partially supported by the
Max Planck Institute in the Mathematical Sciences (Leipzig). Research of S. Nik\v cevi\'c partially supported by MM 1646
(Srbija).

\end{document}